\documentclass[11pt,a4paper,openany,oneside]{amsart}
\usepackage[a4paper,left=3cm,right=3cm, top=3cm, bottom=3cm]{geometry}
\usepackage[latin1]{inputenc}
\usepackage{mathrsfs}
\usepackage[T1]{fontenc}
\usepackage{graphicx}
\usepackage{amsmath}
\usepackage{esint}
\usepackage{mathtools}
\usepackage{amsthm}
\usepackage{amssymb}
\usepackage{hyperref}
\usepackage{color}
\usepackage{enumitem}
\theoremstyle{plain}
\newtheorem{thm}{Theorem}[section]
\newtheorem{cor}[thm]{Corollary}
\newtheorem{lem}[thm]{Lemma}
\newtheorem{prop}[thm]{Proposition}
\theoremstyle{definition}
\newtheorem{defn}{Definition}[section]
\theoremstyle{remark}
\newtheorem{remark}{Remark}

\newcommand{\R}{\mathbb{{R}}}
\newcommand{\N}{\mathbb{{N}}}

\newcommand{\Ne}{\mathcal{N}}
\newcommand{\Holder}{\text{H\"{o}lder}}
\newcommand{\Lp}{L^p}
\newcommand{\Ho}{H_{0}^{1}(\Omega)}
\newcommand{\Omegabar}{\overline{\Omega}}
\newcommand{\dOmega}{{\partial\Omega}}
\newcommand{\dBR}{{\partial B_R}}

 \def\dd{\, {\rm d}}

\newcommand{\iii}[1]{{\left\vert\kern-0.25ex\left\vert\kern-0.25ex\left\vert #1 \right\vert\kern-0.25ex\right\vert\kern-0.25ex\right\vert}}

\usepackage{tikz}
\usepackage{subcaption}
\usepackage{pgfplots}

\pgfplotsset{every axis/.append style={
		axis x line=middle,    
		axis y line=middle,    
		axis line style={<->,color=blue}, 
		xlabel={$x$},          
		ylabel={$y$},          
}}

\begin{document}
	\title[Symmetry and uniqueness for a hinged plate problem in a ball]{Symmetry and uniqueness \\for a hinged plate problem in a ball}
	\author{Giulio Romani}
	\address{Dipartimento di Scienza e Alta Tecnologia Universit\`{a} degli Studi dell'Insubria and RISM-Riemann International School of Mathematics Villa Toeplitz, Via G.B. Vico, 46 - 21100 Varese, Italy}
	\email{giulio.romani@uninsubria.it}

\date{\today}



\begin{abstract}
	In this paper, we address questions related to symmetry, radial monotonicity, and uniqueness for a semilinear fourth-order boundary value problem in the ball of $\R^2$ originating from the Kirchhoff-Love model of deformations of thin plates. We first show the radial monotonicity for a broad class of biharmonic problems. The proof of uniqueness is based on ODE techniques and applies to the whole range of the boundary parameter. For an unbounded subset of this range we also prove symmetry of the ground states by means of a rearrangement argument which makes use of Talenti's comparison principle. This work complements the results in \cite{MIO}, where existence and positivity of solutions were previously analysed.
\end{abstract}

\maketitle

\section{Introduction and main results}

Let
$B_R\subset\R^2$ be the ball of radius
$R>0$ centred at the origin and denote by
$\kappa=\frac1R$ its curvature. 
We aim at investigating qualitative 
properties of solutions of the 
fourth-order semilinear problem
\begin{equation}\label{PDEpSteklov}
	\begin{cases}
		\Delta^2u=|u|^{p-1}u\quad&\mbox{in }B_R,\\
		u=\Delta u-(1-\sigma)\kappa u_n=0
		\quad&\mbox{on }\dBR,
	\end{cases}
\end{equation}
in particular regarding symmetry and 
uniqueness of positive ground state solutions.
Here,
$p\in(0,1)\cup(1,+\infty)$,
$\sigma>-1$ is a boundary parameter, and
$u_n:=\partial_nu$ stands for the 
(outward pointing) normal derivative 
of the function
$u$. Fourth-order and, more generally, 
higher-order problems have attracted 
interest mainly because of the intrinsic 
difficulties that such problems entail. 
First, positivity of the datum in 
general fails to be preserved by the 
solutions, and this property is strictly 
related to the kind of boundary 
conditions (BCs) one has to deal with. 
Consequently, many techniques, 
such as Harnack's inequalities and 
maximum principles, which are familiar 
for second-order equations, do not 
extend to this setting. Another main 
obstacle is that the positive and 
negative part of a function belonging 
to a higher-order Sobolev space do not 
lie within the same space, due to 
possible jumps of the derivatives. 
This implies that also truncations and 
rearrangement techniques cannot be 
directly applied.

In the literature, the equations of order
$2m\geq4$ are in general endowed with 
Dirichlet BCs 
$$
u=\partial_nu=\dots=\partial_n^{m-1}u=0
$$
or Navier BCs 
$$
u=\Delta u=\dots=\Delta^{m-1}u=0 . 
$$
These two cases enjoy opposite features. 
On the one hand, when
$\dOmega$ is sufficiently smooth, with 
Navier BCs one can decouple the problem 
into a system of second-order equations, 
for which the maximum principle holds, 
and therefore the solution inherits the 
sign from the data. On the other hand, 
the positivity preserving is in general 
lost in the Dirichlet case, even for 
smooth and convex domains, except for 
peculiar situations in which one can 
rely on a global analysis of the Green 
function, such as for the case of the 
ball and its smooth deformations. 
For a comprehensive discussion on the 
positivity preserving property for 
higher-order problems we refer to 
\cite{GGS} and for more recent 
developments to \cite{GRS,CTa}.

As a sort of intermediate case 
between Navier (for
$\sigma=1$) and Dirichlet BCs 
(for
$\sigma\to+\infty)$, Steklov BCs 
($u=\Delta u-(1-\sigma)\kappa u_n=0$) 
arise naturally from the model of thin 
hinged plates. The behaviour of a thin 
and hinged plate under the action of 
a vertical external force of density
$f$ can be in fact modeled by the 
Kirchhoff-Love functional
\begin{equation}\label{I}
	I_\sigma(u)=\int_\Omega
	\dfrac{(\Delta u)^2}2-(1-\sigma)
	\int_\Omega\det(\nabla^2u)-\int_\Omega fu,
\end{equation}
where the bounded domain
$\Omega\subset\R^2$ describes the 
shape of the plate and
$u$ its deflection from the original 
unloaded position, see e.g. 
\cite{VK}. 
The parameter
$\sigma$, called Poisson ratio, 
depends on the material and 
measures its transverse expansion 
(resp. contraction), 
according to its positive 
(resp. negative) sign, when subjected 
to an external compressing force. 
In case of a sufficiently smooth boundary
$\dOmega$, the functional
$I_\sigma$ can be rewritten as
\begin{equation}\label{I_kappa}
	I_\sigma(u)=\int_\Omega\dfrac{(\Delta u)^2}
	2-\frac{1-\sigma}2\int_\dOmega\kappa u_n^2
	-\int_\Omega fu,
\end{equation}
see \cite{PS,MIO}, where
$\kappa$ is the signed curvature of the
boundary (that is, positive on convex parts), 
and critical points of
$I_\sigma$ correspond to weak solutions of
\begin{equation}\label{PDE_lin}
	\begin{cases}
		\Delta^2u=f\quad
		&\mbox{in }\Omega,\\
		u=\Delta u-(1-\sigma)\kappa 
		u_n=0\quad&\mbox{on }\dOmega.
	\end{cases}
\end{equation}
Note that the first BC is due to the
fact that the plate is hinged, while
the second BC is a consequence of
the integration by parts, 
see \cite{GS}.

The literature for fourth-order problems 
endowed with Steklov BCs is not wide. 
The linear problem
\eqref{PDE_lin} and 
its positivity preserving property have been 
investigated in \cite{GS,PS}, while its 
semilinear counterpart has been considered in 
\cite{BGM,BGW,GP,MIO}. One finds that the 
issues of existence and nonexistence, as well 
as positivity, are strongly related to the parameter
$\sigma$ and to the properties of the domain. 
Indeed, when
$\Omega\subset\R^2$ is a bounded convex domain of class
$C^{1,1}$ and one considers the semilinear setting
$f(s)=|s|^{p-1}s$ with
$p\in(0,1)\cup(1,+\infty)$, there exist
$\sigma_*\leq-1$ and
$\sigma^*>1$ (possibly infinite) depending on
$\Omega$ such that 
\begin{equation}\label{PDE_NL}
	\begin{cases}
		\Delta^2u=|u|^{p-1}u\quad&\mbox{in }
		\Omega,\\
		u=\Delta u-(1-\sigma)\kappa u_n=0
		\quad&\mbox{on }\dOmega
	\end{cases}
\end{equation}
has no positive solutions for
$\sigma\leq\sigma_*$, while if
$\sigma\in(\sigma_*,\sigma^*)$ it admits a 
positive ground state solution, namely a 
least-energy critical point of the functional
\begin{equation}\label{J_kappa}
	J_\sigma(u)=\int_\Omega\dfrac{(\Delta u)^2}
	2-\frac{1-\sigma}2\int_\dOmega\kappa u_n^2
	-\int_\Omega\frac{|u|^{p+1}}{p+1},
\end{equation}
see \cite[Theorem 1.1]{MIO}. 
We note that for
$\sigma\searrow\sigma_*$ the
$L^\infty$-norm of positive ground 
state solutions blows-up as
$p\in(0,1)$, while their
$H^2$-norm decays to
$0$ when
$p>1$, see \cite[Theorem 1.2]{MIO}. 
When
$\Omega=B_R\subset\R^2$, 
one can refine the analysis and prove in 
addition that
$\sigma_*=-1$ and
$\sigma^*=+\infty$. We point out that the 
latter equality is due to the fact that the 
corresponding Dirichlet problem in the ball 
is positivity preserving. Hence, by means of 
Palais principle of symmetric criticality, 
the existence of a positive radial solution 
can be proved for all
$\sigma>-1$. Moreover, 
it is known that positive radial solutions 
are radially decreasing when
$\sigma\in(-1,1]$, 
\cite[Propositions 7.1 and 7.3]{MIO}. 
However, it was still open whether ground state solutions are radially symmetric, as well as which is the radial behavior of positive radial solutions when $\sigma>1$. It is worth noticing that, differently from the second-order setting, for issues about symmetry and radial behaviour one cannot rely on Gidas-Ni-Nirenberg type results, see \cite{NoGNN}, or direct symmetrisation methods. Indeed, on the one hand, the symmetrised $u^*$ of a function $u\in H^2$ may not belong to the same Sobolev space. On the other hand, for higher-order problems with Dirichlet boundary conditions, the well-established moving-plane method -- widely applicable to second-order problems, even in delicate settings like fractional equations \cite{CLL} or Schr\"odinger-Poisson systems \cite{CW} -- requires highly precise estimates of the Green's function, see \cite{BGW}, which are known only in the case of the ball. We refer to \cite{MR} for an application of the moving-plane method in the context of uniform a priori bounds for higher-order Dirichlet problems.
\vskip0.2truecm
In Section \ref{Section_Decay},
we prove that positive radial solutions of
\eqref{PDEpSteklov} in the whole range
$\sigma>-1$ are decreasing in the radial 
variable -- as a particular case of a more 
general result -- and show that the ground 
state solutions of
\eqref{PDEpSteklov} are radially symmetric 
for an unbounded subset of the range of
$\sigma$. More precisely, we prove the 
following results.
\begin{prop}\label{rad_decr}
	Let
	$N\geq2$,
	$R>0$,
	$f:\R^N\times\R\to\R$ be a smooth 
	positive nonlinearity, and
	$u$ be a smooth positive radial solution of
	\begin{equation}\label{PDE_gen}
		\begin{cases}
			\Delta^2u=f(x,u)\quad&\mbox{in }
			B_R\subset\R^N,\\
			u=0\quad&\mbox{on }\dBR.
		\end{cases}
	\end{equation}
	Then
	$u$ is strictly radially decreasing.
\end{prop}
\begin{thm}\label{symm}
	Let
	$R>0$,
	$p\in(0,1)\cup(1,+\infty)$ and
	$\sigma\geq1$. Then every ground state solution of
	\eqref{PDEpSteklov} is positive, radially symmetric, 
	and decreasing in the radial variable.
\end{thm}

We stress the fact that in Proposition 
\ref{rad_decr} we only require that
$u$ vanishes on the boundary of the ball, 
\textit{without prescribing the second 
	boundary condition}. 
This makes the result applicable for a wide 
class of radially symmetric biharmonic
problems, and in particular to
\eqref{PDEpSteklov} for the whole range
$\sigma>-1$. Moreover, it extends the 
validity of Conjecture 1 in \cite{NoGNN} in case
$f(u)\geq0$ beyond the known cases of 
Navier and Dirichlet boundary conditions, 
see Remarks 1 and 2 therein. 
The positivity part of the statement of 
Theorem \ref{symm} follows directly from 
\cite[Corollary 6.21]{MIO}, while the 
proof of the symmetry is based on Talenti's 
comparison principle 
\cite{Talenti,Kesavan}; 
the restriction
$\sigma\geq1$ is due to the fact that 
the boundary term in
\eqref{I_kappa} changes sign according 
to the sign of
$1-\sigma$, which makes our argument
applicable only in that case. 
We envisage however that symmetry
still holds if
$\sigma\in(-1,1)$.

Another interesting and notoriously 
hard topic concerns uniqueness. For 
the subcritical equation
$\Delta^2u=|u|^{p-1}u$ in balls of
$\R^N$, uniqueness of positive solutions 
has been established by Troy \cite{Troy} 
assuming Navier BCs
\begin{equation}\label{PDEpNav}
	\begin{cases}
		\Delta^2u=|u|^{p-1}u\quad&\mbox{in }B_R\subset\R^N,\\
		u=\Delta u=0\quad&\mbox{on }\dBR,
	\end{cases}
\end{equation}
by decoupling the fourth-order problem as 
a coupled system of second-order 
equations, and by Dalmasso \cite{Dalmasso}
assuming Dirichlet BCs 
\begin{equation}\label{PDEpDir}
	\begin{cases}
		\Delta^2u=|u|^{p-1}u\quad
		&\mbox{in }B_R\subset\R^N,\\
		u=u_n=0\quad&\mbox{on }\dBR,
	\end{cases}
\end{equation}
by means of ODE-techniques. 
The analysis of the latter case 
was then pursued by Ferrero, Gazzola, 
and Weth in \cite{FGW}. 
We summarise here at once their results:
\begin{enumerate}[label=(\textit{\roman*})]
	\item the Navier problem
	\eqref{PDEpNav} has a unique positive 
	solution, which is radially symmetric 
	and radially decreasing;
	\item the Dirichlet problem
	\eqref{PDEpDir} has a unique radial 
	positive solution, which is radially 
	decreasing.  Moreover, the Rayleigh quotient
	$\frac{\|\Delta u\|_2^2}{\|u\|_p^2}$ in
	$H^2_0(B_R)$ has a positive, radial, 
	and radially decreasing minimizer, 
	which is thus 
	(up to a multiplicative constant) unique.
\end{enumerate}
In particular, the technique developed in 
\cite{Dalmasso} seems very ductile, as it 
was employed also to prove uniqueness in 
the context of Lane-Emden systems and 
nonvariational polyharmonic systems, 
see \cite{Dalmasso2,Delia,DD}. We also point out a recent interesting result in \cite{AE}, which proves uniqueness for the subcritical Lane-Emden equation in the ball in the fractional setting, using Morse theory.
\vskip0.2truecm
Our next result concerns the uniqueness 
of positive radial solutions of
\eqref{PDEpSteklov} and applies to the 
whole range for the boundary parameter
$\sigma$.

\begin{thm}\label{Thm_Uniq}
	Let
	$R>0$,
	$p\in(0,1)\cup(1,+\infty)$ and
	$\sigma>-1$. There exists a unique 
	positive radial solution for the 
	Steklov problem
	\eqref{PDEpSteklov}, which is strictly 
	decreasing in the radial variable.
\end{thm}

The proof of Theorem \ref{Thm_Uniq}, 
which can be found in Section \ref{Section_Uniq}, 
follows an ODE-argument based on Dalmasso's 
technique, nontrivially adapted to be 
employed in our setting. This uniqueness 
result may be also regarded as a subcritical 
counterpart of the uniqueness result in 
\cite{GP}, where in dimension
$N\geq5$ the critical nonlinearity
$f(u)=|u|^{2_*-1}u$ with
$2_*:=\frac{2N}{N-4}$ has been considered, 
by means of techniques which are specific 
for the critical case.

As an immediate consequence, in light of 
\cite[Theorem 1.1 and Proposition 7.1]{MIO}, 
and combining Proposition \ref{rad_decr} and 
Theorems \ref{symm} and \ref{Thm_Uniq},
we get the following.
\begin{cor}
	Let
	$R>0$,
	$p\in(0,1)\cup(1,+\infty)$, and
	$\sigma>-1$.
	\begin{itemize}
		\item If
		$\sigma\geq1$ there exists a unique 
		ground state solution of
		\eqref{PDEpSteklov}, which is positive, 
		radially symmetric, and strictly radially
		decreasing.
		\item If
		$\sigma\in(-1,1)$ there exists a unique 
		positive radial solution of
		\eqref{PDEpSteklov}, which is strictly
		radially decreasing.
	\end{itemize}
\end{cor}

This result partially answers the first out of 
three open questions in \cite[Section 8]{MIO}, 
about symmetry, radial behaviour and uniqueness. 
The last question therein concerned the necessity 
of the convexity assumption for the domain in 
proving positivity for ground state solutions of
problem
\eqref{PDE_NL}. For the sake of 
completeness, 
in the last Section \ref{Section_Nonconvex},
we report a result from \cite{MIO_tesi}, 
which shows that the convexity condition is not 
necessary. Indeed, we prove that for a specific 
class of domains which are nonconvex deformations 
of a ball, the model case of which are the 
so-called \textit{lima\c{c}ons de Pascal} defined in 
\eqref{lima\c{c}ons}, 
the ground states of 
\eqref{PDE_NL} 
are positive. We present here the result for the 
lima\c{c}ons, referring to 
Theorem \ref{Thm_nonconvex} 
for a more general statement.

\begin{thm}
	Let
	$p\in(0,1)\cup(1,+\infty)$ and 
	$\Omega_a\subset\R^2$ with
	$a\in[0,\bar a]$ be a 
	lima\c{c}ons de Pascal. 
	Then there exists 
	$\nu_*(\Omega_a)\in[\sigma_*(\Omega_a),1)$ 
	such that the ground state solutions of 
	\eqref{PDE_NL} with
	$\sigma>\nu_*$ are positive in 
	$\Omega_a$.
\end{thm}

\section{Symmetry and radial monotonicity}
\label{Section_Decay}

\subsection{Radial behaviour of positive
	radial solutions}

Let us begin with a simple but important 
well-known property. Henceforth, given a 
radially symmetric function, we denote by
$\cdot'$ its derivative in the radial direction.

\begin{lem}\label{maxprlem}
	Let
	$u:[0,R]\rightarrow\R^+$ be the restriction 
	to the radial variable of a
	$C^2$ radial function defined in
	$B_R\subset\R^N$. Then for
	$t\in[0,R]$ there holds
	\begin{equation}\label{maxpr}
		t^{N-1}u'(t)=\int_0^ts^{N-1}\Delta u(s)\dd s.
	\end{equation}
\end{lem}

\begin{proof}[\bf Proof.]
	It directly comes from integration by parts of 
	the radial representation of the 
	laplacian for radial functions
	\begin{equation*}
		\Delta u(s)=u''(s)+\frac{N-1}su'(s).
		\qedhere
	\end{equation*}
\end{proof}

For convenience, we state here a straightforward 
application of Lemma \ref{maxprlem}. 
From now on, by a little abuse of notation, if
$u$ is a radial function, we still indicate with
$u$ its restriction to the radial variable, i.e.,
$u(r):=u(x)$ for
$r=|x|$.

\begin{lem}\label{der_lapl}
	If
	$u$ is a positive radial solution of
	\eqref{PDE_gen}, then
	\begin{equation}
		(\Delta u)'>0\qquad\mbox{in}\;\;(0,R]
	\end{equation}
\end{lem}

\begin{proof}[\bf Proof.]
	Since
	$f$ is positive, it is sufficient to 
	apply Lemma \ref{maxprlem} with
	$g=\Delta u$.
\end{proof}

\begin{proof}[\bf Proof of Proposition \ref{rad_decr}]
	Suppose by contradiction that there exists
	$a\in(0,R)$ such that
	$u'(a)=0$. Applying Lemma \ref{maxprlem} we get
	\begin{equation*}
		0=u'(a)=\frac1{a^{N-1}}\int_0^as^{N-1}\Delta u(s)\dd s.
	\end{equation*}
	Hence, necessarily
	$\Delta u$ has to change sign in
	$(0,a)$ and, being increasing by Lemma \ref{der_lapl},
	there exists
	$c\in(0,a)$ such that
	$\Delta u<0$ in
	$[0,c)$ and
	$\Delta u>0$ in
	$(c,R]$. Therefore, for any
	$t>a$, we get
	\begin{equation*}
		u'(t)=\frac1{t^{N-1}}\int_0^ts^{N-1}\Delta u(s)\dd 
		s=\frac1{t^{N-1}}\int_a^ts^{N-1}\Delta u(s)\dd s>0.
	\end{equation*}
	Since
	$u$ is positive, this contradicts the
	boundary condition
	$u(R)=0$, since necessarily
	$u'(R)\leq0$. This shows that
	$u'$ does not change sign in
	$(0,R)$, that is,
	$u$ is monotonically decreasing in
	$[0,R]$.
\end{proof}

\subsection{Symmetry of ground state solutions}
Let us first recall the basic definitions 
in the theory of symmetrisation, 
see e.g. \cite{Kesavan}. If
$\Omega\subseteq\R^N$ is a bounded measurable set,
denote by
$\Omega^*$ the ball in
$\R^N$ such that
$|\Omega|=|\Omega^*|$. Moreover, for
$u:\Omega\rightarrow\R$ measurable, 
the Schwarz symmetrisation
$u^*$ of
$u$ is defined as the function
$u^*:\Omega^*\rightarrow\R$ which is radially symmetric, 
non-increasing with respect to the Euclidean norm
$|x|$ for
$x\in\Omega$, and such that
$\big|\{x\in\Omega\,|\,u>t\}^*\big|
=\big|\{x\in\Omega^*\,|\,u^*>t\}\big|$ 
for all
$t\in\R$. The following properties are 
well-known:
\begin{enumerate}
	\item[\textit{i})]
	$u^*$ is radial and radially decreasing;
	\item[\textit{ii})] if
	$u\in\Lp(\Omega)$ with
	$p\geq1$, then
	$u^*\in\Lp(\Omega^*)$,
	$u^*\geq 0$ and
	$\|u\|_{\Lp(\Omega)}=\|u^*\|_{\Lp(\Omega^*)}$.
\end{enumerate}
Notice that
$u\in H^2(\Omega)$ does not imply
$u^*\in H^2(\Omega^*)$. 
This fact prevents a direct application of the 
symmetrisation techniques to the 
higher-order problems. However, one may 
circumvent this drawback by means of 
Talenti's comparison principle which 
we recall in the following form.
\begin{thm}[Talenti's comparison principle, 
	\cite{Talenti,Kesavan}]\label{Talenti}
	Let
	$\Omega\subset\R^N$ be a bounded domain,
	$f\in L^2(\Omega)$, and let
	$u\in H^2(\Omega)\cap\Ho$ and
	$v\in H^2(\Omega^*)\cap H^1_0(\Omega^*)$ 
	be the unique strong solutions of the problems
	\begin{equation}\label{Tal1}
		\left\{\begin{array}{ll}
			-\Delta u=f,\ u\geq0&\mbox{in }\Omega,\\
			u=0&\mbox{on }\dOmega,
		\end{array}\right.
		\quad\mbox{and}\quad
		\left\{\begin{array}{ll}
			-\Delta v=f^*&\mbox{in }\Omega^*,\\
			v=0&\mbox{on }\dOmega^*,
		\end{array}\right.
	\end{equation}
	respectively. Then
	$u^*\leq v$ in
	$\Omega^*$.
	Suppose in addition that
	$\Omega$ is smooth,
	$f\geq 0$, and denote by
	$\mu(t):=|\{u>t\}|=|\{u^*>t\}|$ and
	$\nu(t):=|\{v>t\}|$ the distributional
	functions of
	$u$ and
	$v$ respectively. If
	$u$ is smooth enough and
	$\mu(t)=\nu(t)$ for all
	$t\geq 0$, then
	$\Omega$ is a ball and
	$u$ is radial.
\end{thm}

In the rest of this Section, we prove 
Theorem \ref{symm}, namely that the 
ground-state solutions for the functional
$J_\sigma$
\eqref{J_kappa} associated to problem
\eqref{PDEpSteklov} are radially symmetric. 
We need to distinguish between
$p\in(0,1)$ and
$p>1$, since in the first case we are 
dealing with global minima, while in the 
second case with minima restricted to 
the Nehari manifold.

\begin{proof}[\bf Proof of Theorem \ref{symm}]
	Let
	$\Omega=B_R$ and
	$u$ be a (positive) ground state of
	$J_\sigma$. Denote by
	$u^*$ its symmetrised function and by
	$v$ the solution of
	\begin{equation*}
		\left\{\begin{array}{ll}
			-\Delta v=(-\Delta u)^*&\mbox{in }B_R,\\
			v=0&\mbox{on }\partial B_R.
		\end{array}\right.
	\end{equation*}
	Then we have \begin{equation}\label{stima1}
		\|u\|_{p+1}=\|u^*\|_{p+1}\leq\|v\|_{p+1}
	\end{equation}
	and
	\begin{equation}\label{stima2}
		\|\Delta u\|_2=\|(\Delta u)^*\|_2
		=\|\Delta v\|_2.
	\end{equation}
	Next we compare the boundary terms. 
	Note that in this case
	$\kappa=\frac1R$, thus constant. 
	Using the divergence theorem and the 
	properties of the symmetrised function,
	we get
	$$
	\int_\dBR u_n=\int_{B_R}\Delta u
	=\int_{B_R} \Delta v=\int_\dBR v_n.
	$$
	Since
	$v$ is a radial function, one infers
	$$
	\Big |\int_\dBR v_n  \Big |=|v_n(R)||\dBR|
	=(v_n^2(R)|\dBR|)^\frac12|\dBR|^\frac12
	=  \Big (\int_\dBR v_n^2  \Big )^\frac12|\dBR|^\frac12
	$$
	and, by the
	$\Holder$ inequality,
	$$
	|\dBR|^\frac12  \Big (\int_\dBR u_n^2  \Big )^\frac12
	\geq  \Big |\int_\dBR u_n  \Big |=
	\Big (\int_\dBR v_n^2  \Big )^\frac12|\dBR|^\frac12.
	$$
	For
	$\sigma\geq1$ this yields
	\begin{equation}\label{stima3}
		(1-\sigma)\int_\dBR u_n^2
		\leq(1-\sigma)\int_\dBR v_n^2.
	\end{equation}
	
	Suppose first
	$p\in(0,1)$. The comparison estimates
	\eqref{stima1}-\eqref{stima3} directly imply
	$J_\sigma(u)\geq J_\sigma(v)$, 
	which is a contradiction in case the 
	inequality is strict. Hence,
	$J_\sigma(u)=J_\sigma(v)$, which, 
	combined to the fact that
	$0\leq u^*\leq v$ pointwise, implies
	$u^*=v$. The proof is completed using the 
	second statement in Theorem \ref{Talenti}.
	
	On the other hand, suppose
	$p>1$. Then the ground state
	$u$ of
	$J_\sigma$ is a minimum in the Nehari manifold, 
	namely,
	\begin{equation}\label{Nehari}
		J_\sigma(u)=\min_{\Ne_\sigma}J_\sigma,
	\end{equation}
	where
	$ \Ne_\sigma:=\left\{w\in H^2(B_R)\cap
	H^1_0(B_R)\setminus\{0\}\,\big|\,J_\sigma'(w)w=0\right\}.
	$
	Note that on
	$\Ne_\sigma$ the functional
	$J_\sigma$ reduces to 
	$$
	{J_\sigma}_{|_{\Ne_\sigma}}(w)=
	\Big ( \frac12-\frac1{p+1}
	\Big ) \int_{B_R}|w|^{p+1},
	$$
	and, moreover, for all
	$w\in H^2(\Omega)\cap\Ho\setminus\{0\}$ 
	there exists a unique
	$t^*_w>0$ such that
	$t^*_ww\in\Ne_\sigma$ and
	$J_\sigma(t^*_ww)=\max\left\{J_\sigma(tw)
	\,|\,t>0\right\}$, 
	see \cite[Lemma 3.3]{MIO}.
	
	Suppose by contradiction that
	$u^*\neq v$. In this case, the inequality
	$J_\sigma(u)\geq J_\sigma(v)$ implies that
	$v\notin\Ne_\sigma$. Indeed, if
	$v\in\Ne_\sigma$, then
	$$J_\sigma(v)=
	\Big ( \frac12-\frac1{p+1}
	\Big ) \int_{B_R}|v|^{p+1}>
	\Big ( \frac12-\frac1{p+1}
	\Big ) \int_{B_R}|u|^{p+1}=J_\sigma(u),$$
	which is a contradiction. Hence, there exists
	$0<t^*_v\neq 1$ such that
	$t_v^*v\in\Ne_\sigma$, and so
	\begin{align*} 
		J_\sigma(t^*_vv)&=\frac{(t^*_v)^2}2
		\Big ( \|\Delta v\|_2^2-\frac{1-\sigma}R
		\int_{\partial B_R}v_n^2
		\Big ) -\frac{(t^*_v)^{p+1}}{p+1}
		\|v\|_{p+1}^{p+1}\\
		&<\frac{(t^*_v)^2}2
		\Big ( \|\Delta u\|_2^2-\frac{1-\sigma}
		R\int_{\partial B_R}u_n^2
		\Big ) -\frac{(t^*_v)^{p+1}}{p+1}
		\|u\|_{p+1}^{p+1}\\
		&=J_\sigma(t^*_vu)<J_\sigma(u)=
		\inf_{w\in\Ne_\sigma}J_\sigma(w),
	\end{align*}
	which is again a contradiction. So necessarily
	$u^*=v$ and we can conclude the radial symmetry of
	$u$ as in the previous case.
\end{proof}

\section{Uniqueness}\label{Section_Uniq}
Uniqueness of radial and positive solution 
for problems
\eqref{PDEpNav} and
\eqref{PDEpDir} was obtained by Dalmasso in 
\cite{Dalmasso}. The proof is based on the 
same argument, except for the last step, 
which relies on the different second BC. 
Hence, it seems natural to try to 
extend this result also for problem
\eqref{PDEpSteklov}.

Let us first recall a Gronwall-type 
inequality which will be used in the proof, 
see \cite[Thereom 58]{Dragomir}.
\begin{lem}\label{Dragomir}
	Let the nonnegative function
	$\varphi(t)$ defined on
	$[t_0,+\infty)$ satisfy the inequality
	\begin{equation*}
		\varphi(t)\leq c+
		\int_{t_0}^t\int_{t_0}^sG(t,s,r)\varphi(r)\dd r\dd s,
	\end{equation*}
	where
	$G(t,s,r)$ is a nonnegative
	$C^1$ function for
	$t\geq s\geq r\ge t_0$ and
	$c>0$. Then
	\begin{equation*}
		\varphi(t)\leq c\exp
		\Big \{\int_{t_0}^t
		\Big [\int_{t_0}^sG(s,s,r)\dd r+
		\int_{t_0}^t\int_{t_0}^s
		\frac{\partial G}{\partial s}(s,r,\theta)
		\dd\theta\dd r\Big ]\!\dd s
		\Big \}.
	\end{equation*}
\end{lem}

\begin{proof}[\bf Proof of Theorem \ref{Thm_Uniq}]
	The radial behaviour was proved in 
	Proposition \ref{rad_decr}. 
	Suppose by contradiction there exist 
	two radial positive solutions for
	\eqref{PDEpSteklov}, called
	$u$ and
	$v$, so we can define
	$y,z\in C^\infty([0,R],\R^+)$
	such that
	$y(|x|)=u(x)$,
	$z(|x|)=v(x)$, and
	$y'(0)=z'(0)=0$.
	By Proposition \ref{rad_decr} and 
	Lemma \ref{der_lapl} we know that
	$$
	y'<0,\ \ z'<0\ \,\mbox{in}\ \,(0,R)\quad\ 
	\mbox{and}\quad\ (\Delta y)'>0,\
	\ (\Delta z)'>0\ \,\mbox{in}\ \,(0,R].
	$$
	Let
	$\lambda>0$ be such that
	$\lambda^{4/(p-1)}=y(0)/z(0)$ and define
	$$w(t):=\lambda^{4/(p-1)}z(\lambda t).$$
	By construction
	$w(0)=y(0)$ and
	$w$ satisfies
	\begin{equation}\label{PDEpSteklovw_N=2}
		\begin{cases}
			\Delta^2w=|w|^{p-1}w\quad\mbox{in }
			[0,R/\lambda],\\
			w(R/\lambda)=\Delta w(R/\lambda)-(1-\sigma)
			\lambda \kappa w'(R/\lambda)=0.
		\end{cases}
	\end{equation}
	Define also
	$R(\lambda):=\min\{R,R/\lambda\}$, 
	so that the common interval of 
	definition for
	$y$ and
	$w$ is
	$[0,R(\lambda)]$.  
	The key point of the argument is to show that
	\begin{equation}\label{KEYPOINT_N=2}
		\Delta y(0)=\Delta w(0).
	\end{equation}
	Suppose, indeed, that
	\eqref{KEYPOINT_N=2} holds. 
	Then, by Lemma \ref{maxprlem} and 
	integration by parts, one gets
	\begin{align*} 
		y(t)-w(t)&=\int_0^ts\log  
		\Big (\frac ts  \Big )(\Delta y(s)
		-\Delta w(s))\,\dd s\\
		&=\int_0^ts\log  
		\Big (\frac ts  \Big )\int_0^sr\log
		\left ( \frac sr
		\right ) (y^p(r)- w^p(r))\dd r\dd s.
	\end{align*}
	Let us now distinguish the cases
	$p>1$ and
	$p\in(0,1)$. In the first case, the function
	$x\mapsto x^p$ is locally Lipschitz in
	$\R^+$, so there exists a constant
	$C>0$ such that
	\begin{align}
		\label{gronwalleqD2_N=2} 
		|y(t)-w(t)|&\leq\int_0^ts\log
		\Big  ( \frac ts
		\Big  ) \int_0^sr\log
		\Big  ( \frac sr
		\Big  ) |y^p(r)- w^p(r)|\dd r\dd s \notag \\
		&\leq\,C\int_0^t\int_0^ss\log
		\Big  ( \frac ts
		\Big  ) r\log
		\Big  ( \frac sr
		\Big  ) |y(r)- w(r)|\dd r\dd s\\
		\notag
		&\leq\varepsilon+
		C\int_0^t\int_0^sG(t,s,r)|y(r)- w(r)|\dd r\dd s,
	\end{align}
	for any
	$\varepsilon>0$, where the function
	$$
	G(t,s,r):=s\log\Big(\frac ts\Big)r\log\Big(\frac sr\Big)
	$$
	is positive as
	$0<r<s<t$. Hence, Lemma \ref{Dragomir} yields
	$$
	|y(t)-w(t)|\leq\varepsilon\exp
	\Big  ( C\int_0^t\int_0^t\int_0^s
	\frac{\partial G}{\partial s}
	(s,r,\theta)\dd\theta\dd r\dd s
	\Big  ) .
	$$
	By taking the limit
	$\varepsilon\to0$, one has
	$y(t)=w(t)$ for all
	$t\in[0,R(\lambda)]$. This implies from
	\eqref{PDEpSteklovw_N=2} that
	$\lambda=1$ and in turn that
	$y\equiv z$ on
	$[0,R]$. The uniqueness of the 
	radial positive solution is thus proved.

	The case
	$p\in(0,1)$ can be handled in the same
	way paying attention to the fact that
	$x\mapsto x^p$ is locally Lipschitz in
	$(0,+\infty)$. Hence, since
	$y$ and
	$w$ are positive and decreasing in
	$[0,R(\lambda))$, we can repeat the argument in
	$[0,a]$ for an arbitrary
	$a\in(0,R(\lambda))$, obtaining
	$y\equiv w$ in
	$[0,a]$, and then extend it by continuity in
	$[0,R(\lambda)]$.
	
	The above discussion shows that what 
	is left to prove is the claim
	\eqref{KEYPOINT_N=2}. The first part of 
	the argument works independently of the 
	second boundary condition (Navier, 
	Dirichlet, or Steklov), since it relies only on
	$u_{|\dBR}=0$, and essentially follows 
	\cite{Dalmasso}. Our contribution enters 
	in the second part, where the second 
	boundary condition has a determinant role.
	
	Let us argue again by contradiction and 
	suppose without loss of generality that 
	\begin{equation}\label{KEYPOINT_1_N=2}
		\Delta y(0)<\Delta w(0).
	\end{equation}
	Defining
	$g:=y-w$ we have 
	\begin{equation}\label{g_in_0_N=2}
		g(0)=0\qquad\mbox{and}\qquad\Delta g(0)<0.
	\end{equation}
	First,
	we prove that
	$\Delta g<0$ in
	$[0,R(\lambda)]$. If by contradiction there exists
	$a\in(0,R(\lambda)]$ such that
	$\Delta g<0$ in
	$[0,a)$ and
	$\Delta g(a)=0$, then for
	$t\in(0,a]$ there holds
	\begin{equation}\label{gpos_N=2}
		y(t)-w(t)=g(t)=\int_0^tg'(s)\dd s=
		\int_0^t\frac1s\int_0^rr\Delta g(r)\dd r\dd s<0.
	\end{equation}
	As a result,
	$\Delta^2g=y^p-w^p<0$ in
	$(0,a]$ and
	\begin{align*}
		\Delta g(t)
		=\Delta g(t)-\Delta g(a)
		& =-\int_t^a(\Delta g(s))'s\dd s \\
		&
		=-\int_t^a\frac1s\int_0^rr\Delta^2 g(r)\dd r\dd s>0,
	\end{align*}
	which contradicts
	\eqref{g_in_0_N=2}. We have therefore proved that
	\begin{equation}\label{Lapl_sign_N=2}
		\Delta(y-w)<0\qquad\mbox{on}\quad[0,R(\lambda)].
	\end{equation}
	Now, we want to exclude all possibilities for
	$\lambda$. Firstly, by the first boundary 
	condition, we evaluate
	\begin{equation*}
		(y-w)(R(\lambda))=\begin{cases}
			y(R/\lambda)-w(R/\lambda)=y(R/\lambda)>0 &
			\mbox{if  } \lambda>1,\\
			y(R)-w(R)=0 & \mbox{if  } \lambda=1,\\
			y(R)-w(R)=-w(R)<0 & \mbox{if  } \lambda<1.
		\end{cases}
	\end{equation*}
	This readily implies that
	$\lambda<1$ by
	\eqref{gpos_N=2} with
	$t=R(\lambda)$ and
	\eqref{Lapl_sign_N=2}. Hence,
	$R(\lambda)=\min\{R,R/\lambda\}=R$. 
	In order to find the contradiction 
	also for the case
	$\lambda<1$, we have to rely on the 
	second boundary condition.
	
	\begin{remark}
		\rm
		The conclusion of the argument 
		for the Navier and the Dirichlet problems contained in \cite{Dalmasso} 
		is at this point very easy, just 
		evaluating the second boundary condition in
		$R$. In the \textit{Navier} case
		\eqref{PDEpNav} (i.e.,
		$\sigma=1$), since
		$\Delta y(R)=0$, we get
		\begin{equation*}
			\Delta(y-w)(R(\lambda))=-\Delta w(R)>0,
		\end{equation*}
		which contradicts
		\eqref{Lapl_sign_N=2} and the proof 
		concludes. In the \textit{Dirichlet} case
		\eqref{PDEpDir}, since
		$y'(R)=0$, we get
		\begin{equation*}
			(y-w)'(R(\lambda))=-w'(R)>0,
		\end{equation*}
		as
		$w$ is decreasing in
		$[0,R/\lambda]$. However, this is contrast with
		\begin{equation}
			\label{der_sign_N=2}
			(y-w)'(t)=\frac1t\int_0^ts\Delta(y-w)(s)\dd s<0
		\end{equation}
		again by
		\eqref{Lapl_sign_N=2}. 
		In both cases we have excluded that
		$\Delta y(0)<\Delta w(0)$. Similarly one 
		can also deal with the case
		$\Delta y(0)>\Delta w(0)$. 
		This therefore implies our claim
		\eqref{KEYPOINT_N=2} and thus 
		the proof of Theorem \ref{Thm_Uniq} is 
		completed. However, for the Steklov problem
		\eqref{PDEpSteklov} a more involved 
		analysis has to be performed.
	\end{remark}

	\noindent
	{\bf Continuation of the proof 
		of Theorem \ref{Thm_Uniq}}. 
	We recall that
	$\lambda<1$ and 
	\begin{equation}\label{yStek_N=2}
		\Delta y(R)=(1-\sigma)\kappa  y'(R),
	\end{equation}
	\begin{equation}\label{wStek}
		\Delta w(R/\lambda)=(1-\sigma)\kappa \lambda w'(R/\lambda).
	\end{equation}
	We begin evaluating the boundary condition
	$\Delta \cdot(R)=(1-\sigma)\kappa\cdot'(R)$ in
	$y-w$ and let us first suppose that
	$\sigma\in(-1,1)$. Using the identity
	\eqref{maxpr} one infers
	\begin{align}
		\label{mediaintegrale_N=2} 
		\Delta(y-w)(R) & -(1-\sigma) \kappa (y-w)'(R)\\
		\notag
		&
		=\Delta (y-w)(R)-\frac{(1-\sigma)\kappa }
		R\int_0^Rs\Delta (y-w)(s)\dd s\\
		\notag
		&=\Delta (y-w)(R)-\frac{(1-\sigma)}2\frac1{|B_R|}
		\int_{B_R}\Delta (y-w)(x)\dd x.
	\end{align}
	Moreover,
	$\Delta(y-w)$ is decreasing. Indeed one has
	$$(\Delta(y-w))'(r)=\frac1r\int_0^rs\Delta^2(y-w)(s)\dd s
	=\frac1r\int_0^rs(y^p-w^p)(s)\dd s<0$$
	and
	$\Delta y(0)<\Delta w(0)$ by
	\eqref{KEYPOINT_1_N=2}, therefore,
	\begin{equation}\label{mediaintegrale2_N=2}
		0>\frac1{|B_R|}\int_{B_R}\Delta(y-w)(x)\dd x>\Delta(y-w)(R).
	\end{equation}
	Combining
	\eqref{mediaintegrale_N=2} and
	\eqref{mediaintegrale2_N=2}, 
	\begin{equation*}
		\Delta(y-w)(R)-(1-\sigma)\kappa 
		(y-w)'(R)<\Delta(y-w)(R)  \Big (1-\dfrac{1-\sigma}{2} 
		\Big ).
	\end{equation*}
	As
	$\Delta(y-w)<0$ in
	$[0,R(\lambda)]=[0,R]$ and
	$\sigma\in(-1,1)$, we infer
	\begin{equation}\label{firstestimate_N=2}
		\Delta(y-w)(R)-(1-\sigma)\kappa (y-w)'(R)<0,
		\quad\quad\mbox{for}\;\;\lambda<1.
	\end{equation}
	In the complementary case
	$\sigma>1$, the inequality
	\eqref{firstestimate_N=2} follows simply recalling
	\eqref{Lapl_sign_N=2} and
	\eqref{der_sign_N=2}.

	The delicate task is now to obtain a contradiction with
	\eqref{firstestimate_N=2} in order to exclude also that
	$\lambda\in(0,1)$. In other words, we aim at proving that
	\begin{equation}\label{OBJECT}
		\Delta(y-w)(R)-(1-\sigma)\kappa (y-w)'(R)\geq 0
		\quad\quad\mbox{for }\lambda<1.
	\end{equation}
	By
	\eqref{yStek_N=2}, this is equivalent to prove that
	\begin{equation}\label{puntodipartenza_N=2}
		-\Delta w(R)+(1-\sigma)\kappa w'(R)\geq 0.
	\end{equation}
	Let us first rewrite the problem in a 
	more convenient way on the unit ball. 
	To this aim, define
	$\tilde w(r):=\rho^sw(\rho r)$ with
	$\rho:=\tfrac R\lambda$ and
	$s:=\tfrac4{p-1}$. It is not hard to see by
	\eqref{PDEpSteklovw_N=2} that
	$\tilde w$ satisfies
	\begin{equation}\label{PDEpSteklov_tilde_w_N=2}
		\begin{cases}
			\Delta^2\tilde w=|\tilde w|^{p-1}
			\tilde w\quad\quad\mbox{in }[0,1),\\
			\tilde w(1)=\Delta\tilde w(1)-(1-\sigma)
			\tilde w'(1)=0.
		\end{cases}
	\end{equation}
	Notice also that
	$\tilde w$ enjoys the same properties as
	$w$, meaning
	$\tilde w'(0)=0$,
	$\tilde w$ is decreasing in
	$[0,1]$ and regular, hence
	$\tilde w''(0)\leq 0$. Then, define
	$f:[0,1)\to\R$ so that
	\begin{equation*}
		f(r):=-\tilde w''(r)-\frac\sigma r\tilde w'(r).
	\end{equation*}
	We have
	$$f(\lambda)=\rho^{s+2}
	\left ( -w''(R)-\frac\sigma R\tilde w'(R)
	\right ) =\rho^{s+2}
	\left ( -\Delta w(R)+(1-\sigma)\kappa w'(R)
	\right ) ,$$
	where the last identity comes from the 
	radial representation of the laplacian. 
	Therefore, showing that
	\eqref{puntodipartenza_N=2} holds for any
	$\lambda\in(0,1)$ is equivalent to show that 
	\begin{equation}\label{f_N=2}
		f(r)\geq0\qquad\mbox{for any}\quad r\in(0,1).
	\end{equation}
	Suppose by contradiction that there exists
	$r_0\in(0,1)$ so that
	$f(r_0)<0$. Since
	$f(0)\geq0$ and
	$f(1)=0$ by
	\eqref{PDEpSteklov_tilde_w_N=2}, 
	necessarily there exist
	$0<a<b\leq 1$ so that
	$f(a)=0=f(b)$ and
	$f<0$ in
	$(a,b)$, so that
	$f'(b)\geq0$. Hence,
	\begin{align*} 
		0\leq f'(b)&=-\tilde w'''(b)-\frac\sigma b
		\tilde w''(b)+\frac\sigma{b^2}\tilde w'(b)\\
		&=-(\Delta\tilde w)'(b)+\frac{1-\sigma}b
		\Big  ( \tilde w''(b)-\frac{\tilde w'(b)}b
		\Big  ) \\
		&\!\!\!\!\!\stackrel{f(b)=0}{=}-(\Delta\tilde w)'(b)
		+\frac{1-\sigma}b
		\Big  ( -\sigma\frac{\tilde w'(b)}b
		-\frac{\tilde w'(b)}b
		\Big  ) \\
		&=-\underbrace{(\Delta\tilde w)'(b)}_{>0}+
		\frac{\sigma^2-1}b\underbrace{\tilde w'(b)}_{<0}<0
	\end{align*}
	by Lemma \ref{der_lapl}. For
	$\sigma>1$ this proves
	\eqref{f_N=2} and concludes the proof.
	
	The complementary case
	$\sigma\in(-1,1)$ requires a further analysis. As
	$f$ is smooth, there must be a point
	$s\in(a,b)$ which is a minimum for
	$f$, namely such that
	$f'(s)=0$ and
	$f''(s)\geq 0$. This implies
	\begin{equation}\label{(1)_N=2}
		-\tilde w'''(s)=\frac\sigma s\tilde w''(s)
		-\frac\sigma{s^2}\tilde w'(s)
	\end{equation}
	and
	\begin{equation}\label{(2)_N=2}
		\tilde w^{(iv)}(s)\leq-\frac\sigma s
		\tilde w'''(s)+\frac{2\sigma}{s^2}\tilde w''(s)
		-\frac{2\sigma}{s^3}\tilde w'(s).
	\end{equation}
	Moreover, by Lemma \ref{der_lapl}
	we know that
	$\Delta\tilde w$ is increasing in
	$(0,1)$, i.e.,
	\begin{equation}\label{(2a)_N=2}
		-\tilde w'''(s)-\frac{\tilde w''(s)}s+
		\frac{\tilde w'(s)}{s^2}<0.
	\end{equation}
	From
	\eqref{(1)_N=2} and
	\eqref{(2a)_N=2}, we readily find
	\begin{equation*}
		\frac{\sigma-1}s
		\Big  ( \tilde w''(s)-\frac{\tilde w'(s)}s
		\Big  ) <0,
	\end{equation*}
	which in turns yields
	\begin{equation}\label{(3)_N=2}
		\tilde w''(s)-\frac{\tilde w'(s)}s>0
	\end{equation}
	as
	$\sigma<1$. Then, by
	\eqref{(1)_N=2} and
	\eqref{(2)_N=2}, we get
	\begin{equation*}
		\begin{split}
			0<\tilde w^p(s)&=\Delta^2\tilde w(s)=
			\tilde w^{(iv)}(s)+\frac2s\tilde w'''(s)
			-\frac{\tilde w''(s)}{s^2}+
			\frac{\tilde w'(s)}{s^3}\\
			&\leq\frac1s
			\Big  ( (2-\sigma)\tilde w'''(s)+(2\sigma-1)
			\frac{\tilde w''(s)}s+(1-2\sigma)
			\frac{\tilde w'(s)}{s^2}
			\Big  ) \\
			&=\frac{2-\sigma}{s^2}
			\Big  ( \frac{2\sigma-1}{2-\sigma}-\sigma
			\Big  ) 
			\Big  ( \tilde w''(s)-\frac{\tilde w'(s)}s
			\Big  ) \\
			&=\frac{\sigma^2-1}{s^2}
			\Big  ( \tilde w''(s)-\frac{\tilde w'(s)}s
			\Big  ) ,
		\end{split}
	\end{equation*}
	which is again a contradiction by
	\eqref{(3)_N=2} for
	$\sigma\in(-1,1)$. The proof is then concluded.
\end{proof}

\section{Positivity of ground states for a class of nonconvex domains}
\label{Section_Nonconvex}
The existence of ground state solutions for problem
\eqref{PDE_NL} was established in \cite{MIO} for
$\sigma>\sigma_*(\Omega)$ provided
$\Omega$ is a bounded domain of class
$C^{1,1}$. Such boundary regularity was 
needed to have a well-defined curvature
$\kappa\in L^\infty(\dOmega)$ and thus
to show that the map
\begin{equation}\label{Hsigma_norm}
	u\mapsto\|u\|_{H_\sigma}:=
	\Big  ( \int_\Omega(\Delta u)^2-(1-\sigma)
	\int_\dOmega\kappa u_n^2
	\Big  ) ^\frac12
\end{equation}
defines a norm on
$H^2(\Omega)\cap H^1_0(\Omega)$ 
which is equivalent to the standard one, 
see \cite[Lemma 5.5]{MIO}. The convexity 
assumption on
$\Omega$ comes into play in the proof of 
positivity, since
$\kappa\geq0$ is essential in the 
techniques applied therein, see 
\cite[Propositions 4.9, 4.12, 5.6 and Theorem 6.20]{MIO}. 
However, for
$\sigma=1$,
\eqref{PDE_NL} reduces to a Navier boundary 
value problem, so it can be equivalently 
rewritten as a coupled system of 
second-order equations - as long as minimal 
smoothness assumptions are fulfilled, 
\cite{NS} - and the positivity of the 
ground state solutions is a consequence 
of the iterated maximum principle, 
regardless of the convexity of the domain. 
Moreover, in the limiting Dirichlet case 
(that is, when
$\sigma\to+\infty$), there are examples 
of nonconvex domains for which the 
linear problem is positivity preserving 
(see Lemma \ref{LemmaDAS} below). 
Therefore, in \cite[Section 8]{MIO} 
the question of understanding whether the 
convexity of the domain is a necessary 
assumption for the positivity of ground states of
\eqref{PDE_NL} was posed. In this last 
section, we show that for a class of
nonconvex domains the ground states of
\eqref{PDE_NL} are positive.

For $a\in [0,\tfrac12 ]$, the 
\textit{lima\c{c}ons} of parameter
$a$ is defined as
\begin{equation}\label{lima\c{c}ons}
	\Omega_a:=\left\{(\rho\cos\varphi,\rho\sin\varphi)\in\R^2\,\big|\,0\leq\rho<1+2a\cos\varphi\right\}.	
\end{equation}
\begin{figure}[h!]
	\begin{subfigure}{.17\textwidth}
		\centering
		\begin{tikzpicture}[baseline]
			\begin{axis}[axis equal, hide axis, clip=false,xmin=-2.12,xmax=2.12,ymin=-2.12,ymax=2.12]
				\addplot[domain=0:360,samples=300,color=black,data cs=polar] (x,{0.77*(1 + 0*cos(x))});
			\end{axis}
		\end{tikzpicture}
	\end{subfigure}
	\qquad
	\centering
	\begin{subfigure}{.17\textwidth}
		\begin{tikzpicture}[baseline]
			\begin{axis}[axis equal, hide axis, clip=false,xmin=-2.12,xmax=2.12,ymin=-2.12,ymax=2.12]
				\addplot[domain=0:360,samples=300,color=black,data cs=polar] (x,{0.77*(1 + 0.5*cos(x))});
			\end{axis}
		\end{tikzpicture}
	\end{subfigure}
	\qquad
	\centering
	\begin{subfigure}{.17\textwidth}
		\begin{tikzpicture}[baseline]
			\begin{axis}[axis equal, hide axis, clip=false,xmin=-2.12,xmax=2.12,ymin=-2.12,ymax=2.12]
				\addplot[domain=0:360,samples=300,color=black,data cs=polar] (x,{0.77*(1 + 2*0.40824829046*cos(x))});
			\end{axis}
		\end{tikzpicture}
	\end{subfigure}
	\qquad
	\centering
	\begin{subfigure}{.17\textwidth}
		\begin{tikzpicture}[baseline]
			\begin{axis}[axis equal, hide axis, clip=false,xmin=-2.12,xmax=2.12,ymin=-2.12,ymax=2.12]
				\addplot[domain=0:360,samples=300,color=black,data cs=polar] (x,{0.77*(1 + cos(x))});
			\end{axis}
		\end{tikzpicture}
	\end{subfigure}
	\vskip-1.4truecm\caption{Lima\c{c}ons for resp.\ 
		$a = 0$,
		$\frac14$,
		$\frac{\sqrt6}6$,
		$\frac12\,$.}
\end{figure}
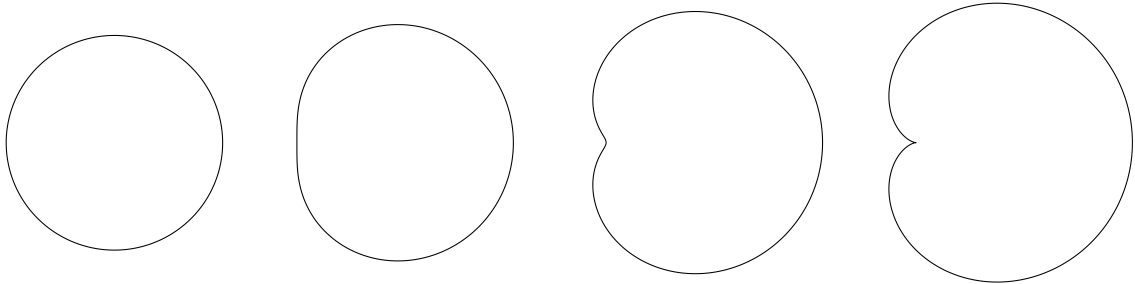
Note that
$\Omega_0$ is the unit disc, and we recall that
$\Omega_a$ is convex if and only if
$a\in\left[0,\tfrac14\right]$. Lima\c{c}ons 
are one of the few examples for which the 
Green function for the Dirichlet problem
\begin{equation}\label{DIRlin}
	\begin{cases}
		\Delta^2u=f\quad&\mbox{in }\Omega,\\
		u=u_n=0\quad&\mbox{on }\dOmega,
	\end{cases}
\end{equation}
is explicit. Moreover, a lower bound of the 
Green function by a positive quantity, which 
is also preserved by small smooth deformations 
of the domain, has been obtained by Dall'Acqua 
and Sweers, see 
\cite[Theorem 5.3.2]{DallAcquaTesi} and 
\cite[Theorem 3.1.3 and Remark 3.4.3]{DASlimacon}.

\begin{defn}
	Let
	$\varepsilon>0$,
	$k\in\N$,
	$\gamma\in(0,1)$. We say that
	$\Omega\subset\R^N$ is
	$\varepsilon$\textit{-close to
		$\Omega'\subset\R^N$ in
		$C^{k,\gamma}$-sense} if there exists a
	$C^{k,\gamma}$ mapping
	$g:\overline{\Omega'}\rightarrow\Omegabar$ such that
	$g(\overline{\Omega'})=\Omegabar$ and
	$\|g-Id\|_{C^{k,\gamma}(\overline{\Omega'})}
	\leq\varepsilon$.
\end{defn}

\begin{lem}[\cite{DallAcquaTesi,DASlimacon}]
	\label{smallperturblimacon}\label{LemmaDAS}
	Let
	$\bar a\in
	( \frac14,\frac{\sqrt6}6
	) $ and
	$\gamma\in(0,1)$. Then there exist
	$\varepsilon_0>0$ and
	$c_1,c_2>0$ such that for every
	$\varepsilon\in[0,\varepsilon_0]$ and
	$a\in[0,\bar a]$ the following holds: if
	$\Omega$ is
	$\varepsilon$-close in
	$C^{2,\gamma}$-sense to the lima\c{c}ons
	$\Omega_a$, then the Green function
	$G_\Omega$ of
	\eqref{DIRlin} satisfies
	\begin{equation}\label{goodestimatelimacons}
		0<c_1D_\Omega(x,y)\leq G_\Omega(x,y)\leq 
		c_2D_\Omega(x,y)\quad\mbox{for all}\
		\,x,y\in\Omega,
	\end{equation}
	where
	$$
	D_{\Omega}(x,y):=d_{\Omega}(x)
	d_{\Omega}(y)\min  \Big \{1,\frac{d_{\Omega}(x)
		d_{\Omega}(y)}{|x-y|^2}  \Big \}\ \mbox{ and } 
	d_{\Omega}(x):=\textrm{dist}(x,\dOmega).
	$$
\end{lem}

Positivity for the \textit{linear} Steklov 
problem has been investigated in a general 
framework by Gazzola and Sweers in \cite{GS}.

\begin{lem}[\cite{GS}, Theorems 4.1 and 2.6]
	\label{GSthm4.1_2.6}
	Let
	$\Omega\subset\R^N$
	$(N\geq 2)$ be a bounded domain of class
	$C^2$ and let
	$0\lneq\beta\in C(\dOmega)$. Then there exist
	$\delta_{1,\beta}=\delta_{1,\beta}(\Omega)\in(0,+\infty)$
	and
	$\delta_{c,\beta}=\delta_{c,\beta}(\Omega)\in[-\infty,0)$
	such that, if
	$\alpha\in C(\dOmega)$ with
	$\delta_{c,\beta}\beta<\alpha\lneq\delta_{1,\beta}\beta$,
	and
	$u\in H^2(\Omega)\cap\Ho$ is a solution of
	\begin{equation}\label{PDEGAZZSteklov}
		\begin{cases}
			\Delta^2u=f\quad&\mbox{in }\Omega,\\
			u=\Delta u-\alpha u_n=0\quad&\mbox{on }\dOmega,
		\end{cases}
	\end{equation}
	then
	$0\lneq f\in L^2(\Omega)$ implies
	$u>0$ in
	$\Omega$. If moreover
	$\dOmega\in C^{4,\gamma}$ for some
	$\gamma\in(0,1)$ and
	\begin{equation}\label{stimaQ}
		G_\Omega(x,y)\geq c\,d_\Omega(x)^2
		d_\Omega(y)^2\quad\mbox{for every}
		\,\,x,y\in\Omega
	\end{equation}
	for some
	$c>0$, where
	$G_\Omega$ is the Green function for
	\eqref{DIRlin}, then
	$\delta_{c,\beta}=-\infty$.
\end{lem}
We recall that the parameter
$\delta_{1,\beta}$ is characterised as 
the first Steklov eigenfunction 
for the problem
\eqref{PDEGAZZSteklov} for
$f=0$ and
$\alpha=\beta$. In other words,
\begin{equation}\label{delta_1beta}
	\delta_{1,\beta}=\delta_{1,\beta}(\Omega):=
	\inf_{u\in H^2(\Omega)\cap H^1_0(\Omega)}
	\frac{\|\Delta u\|_2^2}
	{\int_\dOmega\beta u_n^2}\,.
\end{equation}

Based on this result, we extend the positivity 
result \cite[Theorem 1.1]{MIO} for ground state 
solutions of the semilinear hinged-plate 
problem in convex domains
\eqref{PDE_NL} to a special class of nonconvex 
domains, obtained by small smooth deformations 
of lima\c{c}ons. We consider indeed 
domains of following class:
\begin{enumerate}
	\item[(D)]
	$\Omega\subset\R^2$ is a bounded domain of class
	$C^{4,\alpha}$ which is
	$\varepsilon$-close in
	$C^{2,\gamma}$-sense to a lima\c{c}ons
	$\Omega_a$, with
	$a\in[0,\bar a]$ and
	$\varepsilon\in[0,\varepsilon_0]$, and
	$\bar a$ and
	$\varepsilon_0$ as in Lemma \ref{smallperturblimacon}.
\end{enumerate}
\begin{thm}\label{Thm_nonconvex}
	Let
	$\Omega\subset\R^2$ be a domain 
	satisfying (D) and
	$p\in(0,1)\cup(1,+\infty)$.
	Then there exists
	$\nu_*(\Omega)\in[\sigma_*(\Omega),1)$ 
	such that the ground state solutions of
	\eqref{PDE_NL} with
	$\sigma>\nu_*$ are positive in
	$\Omega$.
\end{thm}

\begin{proof}[\bf Proof.]
	We consider only the case
	$p>1$, since for
	$p\in(0,1)$ the same proof holds with small 
	adaptations. We split the proof in two cases, 
	according to the sign of
	$1-\sigma$.
	
	\smallskip
	
	\noindent
	{\bf Case $\boldsymbol{\sigma\leq1}$.} 
	First, we show that there exists
	$\nu_*$ such that the \textit{linear} problem
	\eqref{PDE_lin} with
	$\sigma>\nu_*$ is positivity preserving. Defining
	$\beta=|\kappa|\in C(\dOmega)$ and
	$\alpha=(1-\sigma)\kappa$, then by 
	Lemma \ref{GSthm4.1_2.6} there exists
	$\delta_{1,|\kappa|}>0$ such that if
	\begin{equation}\label{upperbound|k|}
		(1-\sigma)\kappa(x)\lneq\delta_{1,|\kappa|}
		|\kappa|(x)\quad\mbox{for all}\ \,x\in\dOmega,
	\end{equation}
	then
	\eqref{PDE_lin} is positivity preserving. 
	Note that if
	$x$ belongs to a convex part of the boundary, then
	$\kappa(x)\geq0$ and so from
	\eqref{upperbound|k|},
	we obtain the condition
	$\sigma>1-\delta_{1,|\kappa|}:=\nu_*$; 
	if otherwise
	$\kappa(x)<0$, we need
	$\sigma<1+\delta_{1,|\kappa|}$, 
	which is of course satisfied if
	$\sigma\leq1$.
	
	Moreover, for
	$\sigma\in(\nu_*,1]$, the map
	\eqref{Hsigma_norm} is still an equivalent norm on
	$H^2(\Omega)\cap H^1_0(\Omega)$. By
	\eqref{delta_1beta} one has indeed
	\begin{equation*}
		\|u\|_{H_\sigma}^2\leq\|\Delta u\|_2^2
		+(1-\sigma)\int_\dOmega|\kappa|u_n^2\leq
		\Big  ( 1+\frac{1-\sigma}{\delta_{1,|\kappa|}}
		\Big  ) \|\Delta u\|_2^2,
	\end{equation*}
	and similarly
	\begin{equation*}
		\|u\|_{H_\sigma}^2\geq\|\Delta u\|_2^2
		-(1-\sigma)\int_\dOmega|\kappa|u_n^2\geq
		\Big  ( 1-\dfrac{1-\sigma}{\delta_{1,|\kappa|}}
		\Big  ) \|\Delta u\|_2^2.
	\end{equation*}
	Therefore, denoted by
	$K$ the positive cone in
	$H^2(\Omega)\cap H^1_0(\Omega)$, 
	one can show that its dual cone
	$K^*$ is a subset of the negative cone, 
	and positivity of ground state solutions of
	\eqref{PDE_NL} follows by means of a dual 
	cones decomposition method. For the details,
	we refer to 
	\cite[Propositions 6.18 and 6.19]{MIO}.
	
	\smallskip
	
	\noindent
	{\bf Case
		$\boldsymbol{\sigma>1}$.} The strategy applied so far
	produces an artificial upper bound, namely
	$\sigma<1+\delta_{1,|\kappa|}$. 
	This is unsatisfying since the Dirichlet problem
	\eqref{DIRlin} is positivity preserving 
	in domains for which the condition (D) holds, 
	and thus we expect to retrieve a result 
	devoid of any upper bound for
	$\sigma$. Let us split the functional
	$J_\sigma$ as
	\begin{equation*}\label{JsigmaNEW}
		J_\sigma(u)=\frac12\iii{u}_\sigma^2
		+(1-\sigma)\int_\dOmega\kappa^-u_n^2-
		\int_\Omega\dfrac{|u|^{p+1}}{p+1},
	\end{equation*}
	where
	$\kappa^-:=\max\{0,-\kappa\}=\frac12
	(|\kappa|-\kappa)$ is the negative part 
	of the curvature and
	\begin{equation*}
		\iii{u}_\sigma:=
		\Big  ( \|\Delta u\|_2^2-(1-\sigma)
		\int_\dOmega|\kappa|u_n^2
		\Big  ) ^\frac12.
	\end{equation*}
	It is easy to see that
	$\iii{\cdot}_\sigma$ defines an equivalent norm on
	$H^2(\Omega)\cap\Ho$, since
	$$(u,v)_\sigma:=\int_\Omega\Delta 
	u\Delta v-(1-\sigma)\int_\dOmega|\kappa|u_nv_n$$
	is a scalar product on
	$H^2(\Omega)\cap\Ho$ and moreover, by
	\eqref{delta_1beta}, we have
	$$
	\|\Delta u\|_2^2\leq\iii{u}_\sigma^2
	\leq\|\Delta u\|_2^2+(\sigma-1)\dfrac{\|\Delta u\|_2^2}
	{\delta_{1,|\kappa|}}=  
	\Big [1+\dfrac{\sigma-1}{\delta_{1,|\kappa|}}  \Big ]
	\|\Delta u\|_2^2.
	$$
	As in the former case, first we 
	investigate the linear Steklov problem, 
	this time endowed with the boundary conditions
	$u=\Delta u-(1-\sigma)|\kappa|u_n=0$, 
	since we want to rely on the norm
	$\iii{\cdot}_\sigma$. 
	Applying Lemma \ref{GSthm4.1_2.6} with
	$\alpha=(1-\sigma)|\kappa|$ and
	$\beta=|\kappa|$ it is easy to obtain 
	the positivity preserving property for
	$\sigma\in(1-\delta_{1,|\kappa|},1+
	|\delta_{c,|\kappa|}|)$. 
	Note that for domains in the class (D) 
	one has
	$\delta_{c,|\kappa|}=-\infty$ by 
	Lemma \ref{GSthm4.1_2.6}, and therefore 
	positivity holds in the linear case for
	$\sigma\in(\nu_*,+\infty)$. 
	Again as in \cite[Proposition 6.18]{MIO} 
	one can show that the dual cone 
	$K^*$ of the positive cone
	$K$ in the norm
	$\iii{\cdot}_\sigma$ is formed by strictly 
	negative (or null) functions. 
	
	Suppose now by contradiction that
	$u\in H^2(\Omega)\cap\Ho$ is a 
	sign-changing ground state of
	$J_\sigma$. By dual cones decomposition 
	\cite[Theorem 3.4]{GGS}, one may split
	$u=u_1+u_2$ with
	$u_1\in K$,
	$0\not\equiv u_2\in K^*$ and
	$(u_1,u_2)_\sigma=0$. Defining
	$w:=u_1-u_2>0$, one has
	$$
	w>|u|\quad\mbox{in}\ \Omega,\qquad 
	w_n^2\geq u_n^2\quad\mbox{on}\ 
	\dOmega,\qquad \iii{w}_\sigma^2=\iii{u}_\sigma^2.
	$$
	By \cite[Lemma 3.3]{MIO}, 
	there exists a unique
	$t^*_w>0$ such that
	$t^*_ww\in\Ne_\sigma$, where
	$\Ne_\sigma$ is the Nehari manifold 
	associated to the functional
	$J_\sigma$. Hence,
	\begin{align*}
		J_\sigma(t^*_ww)&=(t^*_w)^2  
		\Big [\frac12\iii{w}_\sigma^2+(\sigma-1)
		\int_\dOmega\kappa^-(-w_n^2)  \Big ]
		-(t^*_w)^{p+1}\int_\Omega\dfrac{|w|^{p+1}}{p+1}\\
		&
		<(t^*_w)^2  \Big [\frac12\iii{u}_\sigma^2
		+(\sigma-1)\int_\dOmega\kappa^-(-u_n^2)  \Big ]
		-(t^*_w)^{p+1}\int_\Omega\dfrac{|u|^{p+1}}{p+1}\\
		&=J_\sigma(t^*_wu)\leq J_\sigma(u),
	\end{align*}
	where the last inequality is due to the fact that
	$u\in\Ne_\sigma$ and thus is the maximum on the half-line
	$\{tu\,|\,t>0\}$. This contradicts the fact that
	$u$ is a ground state solution. Therefore,
	$u=u_1\geq0$. Finally, as
	$u$ is a critical point of
	$J_\sigma$, for each positive test function
	$\varphi\in H^2(\Omega)\cap\Ho$ one has
	$$
	(u,\varphi)_{H_\sigma}=
	\int_\Omega\Delta u\Delta\varphi-(1-\sigma)
	\int_\dOmega\kappa u_n\varphi_n=\int_\Omega u^p\varphi\geq 0,
	$$
	which implies
	$-u\in K^*$. This in particular yields
	$-u<0$, that is,
	$u>0$.
\end{proof}

\noindent
{\bf Acknowledgments.}
The Author wishes to thank Delia Schiera and Enea Parini for many fruitful discussions during the preparation of this paper. The partial support by INdAM-GNAMPA Project 2024 \textit{New perspectives on Choquard equation through PDEs with local sources} (CUP E53C23001670001) is kindly acknowledged.

\end{document}